\documentclass[reqno,12pt]{amsart}

\usepackage[T1]{fontenc}
\usepackage[utf8]{inputenc}
\usepackage{lmodern}
\usepackage{amssymb,mathtools}
\usepackage{mathrsfs}
\usepackage{amscd}
\usepackage{graphicx}
\usepackage{array}
\usepackage{booktabs}
\usepackage{multirow}
\usepackage{xcolor}
\usepackage{microtype}
\usepackage{float}
\usepackage[a4paper,margin=1in]{geometry}
\allowdisplaybreaks[4]
\theoremstyle{plain}
\newtheorem{theorem}{Theorem}[section]

\newtheorem{lemma}[theorem]{Lemma}

\theoremstyle{definition}

\theoremstyle{remark}

\numberwithin{equation}{section}

\newcommand{\Om}{\Omega}

\newcommand{\R}{\mathbb{R}}
\newcommand{\dd}{\,\mathrm{d}}

\newcommand{\C}{\mathbb{C}}

\newcommand{\ii}{\mathrm i}
\newcommand{\supp}{\operatorname{supp}}
\newcommand{\Rea}{\operatorname{Re}}
\newcommand{\Ima}{\operatorname{Im}}
\newcommand{\norm}[1]{\left\|#1\right\|}

\newcommand{\calH}{\mathcal H}
\newcommand{\calY}{\mathcal Y}
\newcommand{\obs}{\mathcal C}
\newcommand{\ip}[2]{\left\langle #1,#2\right\rangle}

\usepackage[
  colorlinks=true,
  linkcolor=blue,
  citecolor=green,
  urlcolor=blue
]{hyperref}

\title[Stability for time-harmonic diffusion]
{Stability for a multi-frequency inverse source problem in diffusion-based tomography}

\author[Qiu]{Dong Qiu}
\address{School of Mathematical Sciences, Zhejiang University,
Hangzhou 310058, China}
\email{qiudong@zju.edu.cn}

\author[Xu]{Xiang Xu}
\address{School of Mathematical Sciences, and Center for
Interdisciplinary Applied Mathematics, Zhejiang University,
Hangzhou 310058, China}
\email{xxu@zju.edu.cn}

\author[R. Zhao]{Ruiting Zhao}
\address{School of Mathematical Sciences, and Center for
Interdisciplinary Applied Mathematics, Zhejiang University,
Hangzhou 310058, China}
\email{zhaoruiting@zju.edu.cn}

\author[Y. Zhao]{Yue Zhao}
\address{School of Mathematics and Statistics, and Key Lab NAA-MOE,
Central China Normal University, Wuhan 430079, China}
\email{zhaoy@ccnu.edu.cn}

\subjclass[2020]{35R30, 35K05, 35J25}
\keywords{inverse source problem, time-harmonic diffusion,
conditional stability, analytic continuation, boundary observability}
\date{}

\begin{document}

\begin{abstract}
This paper is concerned with the inverse source problem for the time-harmonic diffusion equation arising in optical tomography. In inhomogeneous media, we establish a conditional logarithmic stability estimate by using multi-frequency boundary flux measurements. The proof combines a boundary observability inequality for the associated heat equation with quantitative analytic continuation. For constant coefficients, we further derive an increasing stability estimate with explicit dependence on the frequency bandwidth. In this case, the analysis employs complex polynomial interpolation of the Fourier–Laplace transform of the source. We also implement a multi-frequency Landweber iteration with a discrepancy-based stopping rule to reconstruct the sources. Numerical experiments show an overall improvement in reconstruction accuracy and resolution as the frequency bandwidth increases under the tested noise settings.
\end{abstract}

\maketitle

\section{Introduction}

Optical tomography is an important imaging modality in medical imaging, where visible and near infrared
light is used to transilluminate tissue, and the resulting measurements of intensity on the tissue
boundary are used to reconstruct the optical properties within the tissue. The widely used reconstruction method is based
on recovering the coefficients of a partial differential equation from boundary measurements that is either the radiative transport equation or its diffusion approximation, which models the spatial
density of photons \cite{AL, Bal2, Bal, Bal1}.

In this paper, we consider the diffusion approximation to the radiative transfer
equation, which in the frequency domain is modeled by
\begin{align}\label{eq:bvp}
-\nabla\!\cdot\!\bigl(D(x)\nabla u(x,\omega)\bigr)
 +\bigl(\mu(x)+\ii\omega\rho(x)\bigr)u(x,\omega)
 =f(x),\qquad x\in\Om,
\end{align}
with Dirichlet boundary condition $u=0$ on $\partial\Omega$. Here $\omega>0$ is the frequency, $f$ is the source function, and $\Om\subset\R^d$, $d\ge2$ is a bounded domain with $C^\infty$ boundary. 
We assume that the coefficients $D,\rho,\mu\in C^\infty(\overline\Om)$, characterizing the optical properties
within the tissue, are positive real-valued functions with bounds
\begin{equation}\label{eq:elliptic}
 0<D_-\le D\le D_+,
 \qquad 0<\rho_-\le\rho\le\rho_+,
 \qquad 0<\mu_-\le\mu\le\mu_+.
\end{equation}
The nonuniqueness result of Arridge and Lionheart concerns simultaneous recovery of the diffusion and absorption coefficients from steady-state boundary data ~\cite{AL}.  Nonuniqueness can also occur in frequency-domain and time-domain optical tomography when the refractive index is unknown. The role of multi-frequency data therefore depends on which coefficients are assumed known.

Assuming the coefficients in the governing equation are known, we study the inverse problem of determining the source from the boundary measurements $D\partial_\nu u\big|_{\partial\Om}$. It should be noted that due to the gauge invariance, one may not uniquely recover all the coefficients and sources even if multi-frequency data are used.
Inverse source theory for the wave equations has been well-established, and a well-known fact is that uniqueness cannot be guaranteed for such problems at a single frequency. On the other hand, the use of multiple frequencies can regain the uniqueness and improve the stability in the reconstruction processes \cite{BLZ, IL, LZZ}. 

The aim of this paper is to establish uniqueness and stability for the inverse source problem for the diffusion equation in the frequency domain.
As with wave equations, uniqueness does not hold at a single
frequency. A typical example is that if $f=L\chi\neq 0$, where $L$ is the PDE operator and $\chi\in C_0^\infty(\Omega)$, then no data will be detected on the boundary,
which implies the non-uniqueness. Thus, multi-frequency data shall be utilized. Since the diffusion equations exhibit properties which differ from the wave equation dramatically,
the methodologies developed for wave equations such as using the plane waves to obtain integral equations (cf. \cite{BLZ, LZZ}) are no longer applicable here. It is worth mentioning that
inverse source theory for diffusion equations in the time domain (parabolic equations) has been well established; see~\cite[Chapter~6]{Isakov} for general background and~\cite{Yamamoto2009} for Carleman-estimate methods.

Related inverse problems for diffusion models have been studied from several perspectives. For coefficient identification, frequency-dependent stability bounds have been derived from Dirichlet-to-Neumann maps under suitable a priori assumptions~\cite{Lai2014}. In time-harmonic optical tomography, Lipschitz stability of the boundary absorption coefficient and H\"older stability of its derivatives have been proved when scattering is known~\cite{DGLN2022,CGNS2023}. In the time domain, logarithmic stability has been obtained for recovering the spatial factor of a heat source from a single partial boundary measurement on an infinite waveguide~\cite{KSS2020}.

For the diffusion equation in the frequency domain with variable coefficients, we prove a conditional logarithmic-type stability estimate for the inverse source problem using multi-frequency
Neumann data, assuming the sources belong to the a priori function space. The analysis combines an observability inequality for the heat equation and a quantitative
analytic continuation result linking data at small and large frequencies. Further, for the diffusion equation with constant coefficients we establish an increasing stability,
which shows that as the upper bound of the frequency band or the regularity of the sources increases, the stability will improve. In this case, we employ the Fourier-Laplace
transform of the source function and construct its complex-valued interpolation polynomial. Numerical experiments are presented to reconstruct the source, which show
that better resolution can be achieved as the frequency increases and thus confirm
the theoretical increasing stability.

This paper is organized as follows. In Section \ref{dp} we establish the well-posedness of the direct problem. 
Section \ref{ip} is devoted to stability estimates for the inverse source problem. Numerical experiments are presented in Section \ref{sec:numerics} to reconstruct the sources
by using the Landweber iteration.

\section{Direct problem}\label{dp}

In this section, we establish the well-posedness of the direct problem and the holomorphicity of the data used for the inverse source problem.

Introduce some functional spaces
\[
 \calH=L^2(\Om,\rho\dd x),\qquad \calY=L^2(\partial\Om),
 \qquad \mathcal D(B)=H^2(\Om)\cap H_0^1(\Om),
\]
and let $B=\rho^{-1}(-\nabla\!\cdot D\nabla+\mu)$.
Thus, we have that $B$ is the strictly positive self-adjoint operator in $\calH$ associated with the variational  form
\[
 \mathfrak a(v,w)=\int_\Om
 \bigl(D\nabla v\cdot\nabla\overline w+\mu v\overline w\bigr)\dd x,
 \qquad v,w\in H_0^1(\Om).
\]
Let $\{\phi_j\}_{j\ge1}$ be an orthonormal eigenbasis in $\calH$ and let $0<\lambda_1\le\lambda_2\le\cdots$ be the eigenvalues, counted with multiplicity. Put
\begin{equation}\label{eq:gdef}
 g=\rho^{-1}f,\qquad
 \obs v=D\partial_\nu v\big|_{\partial\Om}.
\end{equation}
Then the equation \eqref{eq:bvp} can be written as $(B+\ii\omega)u=g$.

\begin{theorem}\label{thm:direct}
For every $g\in\calH$ and $\omega\in\R$, \eqref{eq:bvp} has a unique solution $u(\omega)\in\mathcal D(B)$, given by
\begin{equation}\label{eq:directres}
 u(\omega)=(B+\ii\omega)^{-1}g
 =\sum_{j=1}^\infty\frac{g_j}{\lambda_j+\ii\omega}\phi_j,
 \qquad g_j=\ip{g}{\phi_j}_{\calH}.
\end{equation}
The co-normal trace $\obs u(\omega)$ belongs to $\calY$.
\end{theorem}
\begin{proof}
Due to the existence of the imaginary part in the equation, we have $|\lambda_j/(\lambda_j+\ii\omega)|\le1$, which yields
the convergence of the sequence \eqref{eq:directres} in $\mathcal D(B)$. From the strict positivity of $B$, we have that the boundary value problem admits 
a unique solution $u(\omega)\in\mathcal D(B)$, and the co-normal trace result follows from the elliptic regularity. The proof is complete.
\end{proof}

\begin{lemma}\label{lem:resolvent}
Let $\beta\in(3/4,1)$. The operator $\obs B^{-\beta}:\calH\to\calY$ is bounded. Letting
\begin{equation}\label{eq:Fdef}
 F(z)=\obs(B+z)^{-1}g,\qquad \Rea z\ge0,
\end{equation}
one has that $F$ is holomorphic in the right-half complex plane and extends holomorphically through the imaginary axis. Moreover, there is a constant $C_\beta$ depending on $\beta$, $\Om$, and the fixed coefficients $D,\rho,\mu$ such that the following estimates hold for $\Rea z\ge0$:
\begin{align}
 \norm{F(z)}_{\calY}&\le C_\beta\norm g_{\calH}(1+|z|)^{\beta-1},\label{eq:Fbound}\\
 \norm{F'(z)}_{\calY}&\le C_\beta\norm g_{\calH}(1+|z|)^{\beta-2}.
 \label{eq:Fprime}
\end{align}
\end{lemma}
\begin{proof}
The spectral theorem gives $\mathcal D(B^\beta)=[\calH,\mathcal D(B)]_\beta$ with equivalent norms. Here $[X_0,X_1]_\beta$ denotes the complex interpolation space between $X_0$ and $X_1$ with parameter $\beta$. Elliptic regularity identifies the graph norm on $\mathcal D(B)$ with the $H^2$ norm, and interpolation therefore gives $\mathcal D(B^\beta)\hookrightarrow H^{2\beta}(\Om)$. Since $2\beta>3/2$, the normal trace maps this space continuously into $H^{2\beta-3/2}(\partial\Om)\hookrightarrow\calY$. The boundedness of $\obs B^{-\beta}$ follows from this trace estimate and the strict positivity of $B$; see~\cite{LM} for interpolation and trace theory.

For $\lambda\ge\lambda_1$ and $\Rea z\ge0$, one has $|\lambda+z|\ge c(\lambda+|z|)$. Consequently, for $m=1,2$, we have
\[
 \sup_{\lambda\ge\lambda_1}\frac{\lambda^\beta}{|\lambda+z|^m}
 \le C(1+|z|)^{\beta-m}.
\]
Then by factoring $\obs(B+z)^{-m}$ as $\obs B^{-\beta}$ times $B^\beta(B+z)^{-m}$ and applying the spectral theorem we prove \eqref{eq:Fbound}. 
Moreover, the estimate \eqref{eq:Fprime} is a consequence of $F'(z)=-\obs(B+z)^{-2}g$.
Finally, as the poles of $(B+z)^{-1}$ lie at $-\lambda_j$, strictly to the left of the imaginary axis, we have the holomorphic continuation through the imaginary axis.
\end{proof}

\section{Inverse source problem}\label{ip}

In this section, we study the inverse source problems. First, we establish a logarithmic stability for the inverse source problem.
Second, assuming the coefficients in the diffusion equation are constants, we prove an increasing stability, which shows that the stability will improve
as the bandwidth of the frequency increases.  

\subsection{An observation inequality for the heat equation}

Consider the time-dependent heat equation  with Dirichlet boundary condition
\begin{equation}\label{eq:heat}
 \rho\partial_t w-\nabla\!\cdot(D\nabla w)+\mu w=0,
 \qquad w|_{\partial\Om}=0.
\end{equation}
The initial condition is $w(0)=g:=\rho^{-1}f$, and thus the solution is
$w(t)=e^{-tB}g$~\cite{Pazy}. For $t>0$, the spectral theorem and Lemma~\ref{lem:resolvent} give the estimate
\begin{equation}\label{eq:heat_trace}
 \norm{\obs (e^{-tB}g)}_{\calY}
 \le C\norm g_{\calH}t^{-\beta}e^{-\lambda_1t/2}.
\end{equation}
Since $\beta<1$, the right-hand side is integrable near zero. Thus the following integral is well-defined in $\calY$ for  $\Rea z>0$:
\begin{equation}\label{eq:Laplace}
 F(z)=\int_0^\infty e^{-zt}\obs (e^{-tB}g)\dd t.
\end{equation}
The identity holds first for finite eigenfunction sums. Passing to $g\in\calH$ is justified on the left by Lemma~\ref{lem:resolvent} and on the right by~\eqref{eq:heat_trace} and dominated convergence in $\calY$.

The following parabolic boundary observability inequality is useful in the subsequent analysis. 

\begin{lemma}\label{lem:heatobs}
There is a constant $C_{\rm obs}$ such that the following observability inequality  holds:
\begin{equation}\label{eq:heatobs}
 \norm{e^{-2B}g}_{\calH}^2
 \le C_{\rm obs}\int_1^2
 \norm{\obs e^{-tB}g}_{\calY}^2\dd t,
 \qquad g\in\calH.
\end{equation}
\end{lemma}
\begin{proof}
The analysis employs the boundary Carleman estimate in \cite[Lemma 1.1]{Im}. Dividing by $\rho$, the heat operator becomes
\[
 \partial_t-\nabla\!\cdot(\alpha\nabla)+b_0\cdot\nabla+c_0,
\]
where
\[
\alpha=D/\rho,\quad b_0=-D\rho^{-2}\nabla\rho,
 \quad c_0=\mu/\rho.
\]
For $0<t<1$ we choose $\psi(x)=x_1+C_\psi>0$ on $\overline\Om$. Let $m=\max_{\overline\Om}\psi$, and set
\[
 \varphi(t,x)=\frac{e^{\lambda\psi(x)}}{t^2(1-t)^2},
 \qquad
 \Phi(t,x)=\frac{e^{\lambda\psi(x)}-e^{2\lambda m}}{t^2(1-t)^2}<0.
\]
Notice  that the function $\psi$ has no critical points. Let $\lambda,\sigma$ be the Carleman parameters,  which are sufficiently large. For a smooth solution $v$ of the homogeneous heat equation with zero boundary value, using the Carleman estimate in \cite[Lemma 1.1]{Im} we have
\begin{equation}\label{eq:carlemanshort}
 \int_0^1\!\int_\Om\varphi^3e^{\sigma\Phi}|v|^2\dd x\dd t
 \le C\int_0^1\!\int_{\partial\Om}
             \varphi e^{\sigma\Phi}|\partial_\nu v|^2\dd S\dd t.
\end{equation}
The interior weight is bounded below on $\overline\Om\times[1/3,2/3]$, and the boundary weight is bounded above on $\partial\Om\times(0,1)$. Then equivalence of the weighted norms and the positive bounds for $D$ yield
\[
 \int_{1/3}^{2/3}\norm{v(t)}_{\calH}^2\dd t
 \le C\int_0^1\norm{\obs v(t)}_{\calY}^2\dd t.
\]
Since $\partial_t\norm{v(t)}_{\calH}^2=-2\mathfrak a(v(t),v(t))\le0$, the left side is at least $\norm{v(1)}_{\calH}^2/3$. The proof of ~\eqref{eq:heatobs} follows by letting $v(t)=e^{-(1+t)B}g$. This solution is smooth up to the boundary on $0\le t\le1$ by parabolic smoothing. For complex-valued $g$, we apply the real-coefficient Carleman estimate separately to the real and imaginary parts and add the inequalities to obtain ~\eqref{eq:heatobs}. The proof is complete.
\end{proof}

\subsection{Stability for the inverse source problem}

Fix $K_0>0$ and assume throughout this subsection that $K\ge K_0$. The domain and coefficients are fixed. Constants in this subsection may depend on $\Om,D,\rho,\mu,s$, and $K_0$, but not on $K,M$, or the data discrepancy.

Introduce the multi-frequency data discrepancy in the finite interval $(0, K)$ by
\begin{equation}\label{eq:delta}
 \delta_K^2=\int_0^K\norm{\obs u(\omega)}_{\calY}^2\dd\omega.
\end{equation}
For $s>0$ and $M\ge0$, we also assume sources belong to the following a priori class
\begin{equation}\label{eq:prior}
\mathcal A_M =\{ g\in\mathcal D\bigl((I+B)^{s/2}\bigr): \|(I+B)^{s/2}g\|_{\calH}\le M\},
\end{equation}

The following lemma provides the continuation of the given data to those at large frequencies. 
\begin{lemma}\label{lem:continue}
Let $g\in\mathcal A_M$ with $M>0$, and let $F$ be given by~\eqref{eq:Fdef}. There are constants $C_0,C_1,c_0>0$, depending only on the fixed background and $K_0$, such that, by setting $N=C_0M$ and $\eta=C_1\delta_K$, if $0<\eta<N$, then for every $R\ge1$ one has
\begin{equation}\label{eq:continue}
 \sup_{|\tau|\le R}\norm{F(1+\ii\tau)}_{\calY}
 \le N\exp\left(-\frac{c_0L}{1+R^2}\right),
\end{equation}
where $L=\log(N/\eta)$.
\end{lemma}
\begin{proof}
By Lemma~\ref{lem:resolvent}, we assume $C_0$ is sufficiently large such that $\norm{F(z)}_{\calY}\le N$ for $\rm{Re} z\geq 0$. Using Chebyshev's inequality we have a measurable set $E\subset(0,K_0)$ with $|E|\ge K_0/2$ such that 
\[
 \norm{F(\ii\omega)}_{\calY}
 \le (2/K_0)^{1/2}\delta_K=\eta, \quad \omega\in E.
\]
Then we let $C_1=(2/K_0)^{1/2}$. For any $y\in\calY$ with $\norm y=1$, the scalar function $h_y(z)=\ip{F(z)}y_{\calY}$ is bounded and holomorphic. The two-constants principle in the right half-plane~\cite{Ransford} yields
\begin{align}\label{tcp}
 |h_y(a+\ii\tau)|\le
 N^{1-\vartheta_E(a,\tau)}\eta^{\vartheta_E(a,\tau)},
\end{align}
where
\[
\vartheta_E(a,\tau)=\frac1\pi\int_E
 \frac{a}{a^2+(\tau-\omega)^2}\dd\omega.
\]
Specifically, the estimate \eqref{tcp} follows by applying the maximum principle to the subharmonic function $\log(|h_y|/N)$, and using the harmonic extension of the boundary majorant $\log(\eta/N)\mathbf1_E$, where $\mathbf1_E$ is the characteristic function of the set $E$.

At $a=1$ and $|\tau|\le R$, we have
\[
 \vartheta_E(1,\tau)
 \ge\frac{K_0}{2\pi(1+(R+K_0)^2)}
 \ge\frac{c_0}{1+R^2}.
\]
This yields~\eqref{eq:continue} according to the definition of $h_y$. The proof is complete.
\end{proof}

\begin{lemma}\label{lem:invert}
Under the assumptions of Lemma~\ref{lem:continue}, let $\beta=7/8$ and $L=\log(N/\eta)$. For all sufficiently large $L$, we have
\begin{equation}\label{eq:heatsmall}
 \sup_{1\le t\le2}\norm{\obs e^{-tB}g}_{\calY}
 \le CM L^{-1/32}.
\end{equation}
\end{lemma}
\begin{proof}
We first justify the inverse Laplace identity (see~\cite{Pazy} for the standard semigroup setting)
\begin{equation}\label{eq:bromwich}
 \obs e^{-tB}g
 =\frac{e^t}{2\pi}\lim_{R\to\infty}
 \int_{-R}^R e^{\ii\tau t}F(1+\ii\tau)\dd\tau,
 \qquad t>0.
\end{equation}
In view of Lemma~\ref{lem:resolvent}, for $|\tau|\ge1$ one has
\begin{align}\label{est}
 \norm{F(1+\ii\tau)}_{\calY}\le CM|\tau|^{\beta-1},
 \qquad
 \norm{\partial_\tau F(1+\ii\tau)}_{\calY}
 \le CM|\tau|^{\beta-2}.
\end{align}
For $S>R\ge1$, integration by parts on the finite interval gives
\[
\begin{split}
 \int_R^S e^{\ii\tau t}F(1+\ii\tau)\dd\tau
 &=\left[\frac{e^{\ii\tau t}F(1+\ii\tau)}{\ii t}\right]_R^S-\frac1{\ii t}\int_R^S
       e^{\ii\tau t}\partial_\tau F(1+\ii\tau)\dd\tau.
\end{split}
\]
Since $\beta<1$, the endpoint at $S$ tends to zero and the derivative integral converges absolutely in $\calY$. Thus
\begin{equation}\label{eq:tailLaplace}
 \left\|\int_R^\infty e^{\ii\tau t}F(1+\ii\tau)\dd\tau\right\|_{\calY}
 \le \frac{CM}{t}R^{\beta-1}.
\end{equation}
The same bound holds for the negative tail. Together with continuity on compact intervals, this proves convergence of the improper integral in~\eqref{eq:bromwich}.

To identify its value, set $g_n=\sum_{j=1}^n g_j\phi_j$ and $F_n(z)=\obs(B+z)^{-1}g_n$. Scalar Laplace inversion of $(\lambda_j+z)^{-1}$ proves~\eqref{eq:bromwich} for each finite sum $g_n$. Let $h_n=g-g_n$. Lemma~\ref{lem:resolvent} and the preceding integration-by-parts argument apply to $h_n$ with $M$ replaced by $\norm{h_n}_{\calH}$. Splitting at $|\tau|=1$ therefore gives
\[
 \left\|\lim_{R\to\infty}\int_{-R}^R e^{\ii\tau t}
       (F-F_n)(1+\ii\tau)\dd\tau\right\|_{\calY}
 \le C(1+t^{-1})\norm{h_n}_{\calH}.
\]
Also,~\eqref{eq:heat_trace} gives $\norm{\obs e^{-tB}h_n}_{\calY}\le C_t\norm{h_n}_{\calH}$. Passing to $n\to\infty$ establishes~\eqref{eq:bromwich}, with convergence uniform for $1\le t\le2$ in these estimates.

For $1\le t\le2$, the correct splitting is
\[
 \obs e^{-tB}g=\frac{e^t}{2\pi}
 \left(\int_{-R}^R+\int_{|\tau|>R}\right)
       e^{\ii\tau t}F(1+\ii\tau)\dd\tau.
\]
Using~\eqref{eq:continue} and~\eqref{eq:tailLaplace} yields
\begin{equation}\label{eq:heatR}
 \norm{\obs e^{-tB}g}_{\calY}
 \le CM\left(
 R\exp\left(-\frac{c_0L}{1+R^2}\right)+R^{\beta-1}
 \right).
\end{equation}
Choose $R=L^{1/4}\ge1$. The first term is bounded by $CL^{1/4}e^{-c\sqrt L}$ and the second term is equal to $L^{-(1-\beta)/4}=L^{-1/32}$, which gives~\eqref{eq:heatsmall}.
The proof is complete.
\end{proof}

In the following theorem, we present the stability estimate for the inverse source problem, which implies the uniqueness of the inverse problem.
\begin{theorem}\label{thm:variable}
Assume $K\ge K_0$ and $g=\rho^{-1}f\in\mathcal A_M$.
We have the following stability estimate for the inverse source problem:
\begin{equation}\label{eq:mainvar}
 \norm f_{L^2(\Omega)}
 \le CM\left(\log\!\left(2+\log\!\left(2+\frac M{\delta_K}\right)\right)\right)^{-s/2}.
\end{equation}
The constant $C$ is independent of $K,M$, and $\delta_K$. \end{theorem}
\begin{proof}
 Fix $L_*>1$ large enough for Lemma~\ref{lem:invert} and the logarithmic comparison below, and choose $\epsilon_*>0$ so small that $\log(C_0/(C_1\epsilon_*))\ge L_*$. If $\delta_K/M\ge\epsilon_*$, then
\[
 \log\!\left(2+\log\!\left(2+\frac M{\delta_K}\right)\right)
 \le \log\!\left(2+\log(2+\epsilon_*^{-1})\right).
\]
The bound $\norm f_{L^2}\le\rho_+^{1/2}M$ proves~\eqref{eq:mainvar} in this case after enlarging $C$. It remains to consider $0<\delta_K/M<\epsilon_*$. In this case, $L=\log(C_0M/(C_1\delta_K))\ge L_*$ is a parameter determined by the relative error. Combining Lemmas~\ref{lem:heatobs} and~\ref{lem:invert} gives
\begin{equation}\label{eq:finalsmall}
 \norm{e^{-2B}g}_{\calH}\le CM L^{-1/32}.
\end{equation}
Then letting $P_\Lambda$ be the spectral projection of $B$ onto $\{\lambda_j\le\Lambda\}$, which yields
\[
 \norm{P_\Lambda g}_{\calH}
 \le e^{2\Lambda}\norm{e^{-2B}g}_{\calH},
\]
and the assumption $g\in\mathcal A_M$ gives
\begin{equation}\label{eq:spectraltail}
 \norm{(I-P_\Lambda)g}_{\calH}^2
 =\sum_{\lambda_j>\Lambda}|g_j|^2
 \le M^2(1+\Lambda)^{-s}.
\end{equation}
Thus, we have
\[
 \norm g_{\calH}\le CM e^{2\Lambda}L^{-1/32}
       +M(1+\Lambda)^{-s/2}.
\]
Letting $\Lambda=(\log L)/128$ we obtain
\[
 \norm g_{\calH}\le CM(\log L)^{-s/2}.
\]
For sufficiently small $\delta_K/M$, this is equivalent to~\eqref{eq:mainvar} up to some constants noting $\norm f_{L^2}\le\rho_+^{1/2}\norm g_{\calH}$.
The proof is complete.
\end{proof}

\subsection{Increasing stability for the inverse source problem with constant coefficients}

In this subsection, we let $D,\mu,\rho$ be positive constants and then consider the diffusion equation
\begin{equation}\label{eq:constant}
 (-\Delta+a+\ii b\omega)u=p\quad\text{in }\Om,
\end{equation}
where
\[
a=\mu/D,\quad b=\rho/D,\quad p=f/D.
\]
For $s>0$ and $M\ge0$, assume the source belongs to the a priori class
\begin{equation}\label{eq:fourierprior}
\mathcal B_M=\{ p\in H^s(\R^d),\quad \supp p\subset\overline\Om,
 \quad\norm p_{H^s(\R^d)}\le M\}.
\end{equation}
Here $p$ denotes the zero extension of the source to $\R^d$, and the $H^s$ norm is the whole-space norm. Define
\[
 d(\omega)=\norm{\partial_\nu u(\omega)}_{L^2(\partial\Om)}.
\]
Define the data discrepancy by
\begin{equation}\label{eq:epsilon}
\varepsilon_K^2=\int_0^K \norm{\partial_\nu u(\omega)}_{L^2(\partial\Om)}^2\dd\omega.
\end{equation}

Choose $R_\Om$ with $\Om\subset B_{R_\Om}(0)$ and introduce the following Fourier-Laplace transform of the source function for $\theta\in S^{d-1}$:
\begin{equation}\label{eq:entireF}
 P_\theta(z)=\int_\Om p(x)e^{-\ii z\theta\cdot x}\dd x,
 \qquad z\in\C.
\end{equation}
The complex function $ P_\theta(z)$ is entire and can be estimated by
\begin{equation}\label{eq:entiregrowth}
 |P_\theta(z)|\le |\Om|^{1/2}M e^{R_\Om|\Ima z|}.
\end{equation}
Let $z(\omega)=\sqrt{-a-\ii b\omega}$ be chosen with $\Rea z(\omega)>0$ and $\Ima z(\omega)<0$ for $\omega>0$. The plane wave $v=e^{-\ii z(\omega)\theta\cdot x}$ with complex frequency satisfies $(-\Delta+a+\ii b\omega)v=0$. Multiplying both sides of the equation by the plane wave and integrating by parts, we obtain the integral identity
\begin{equation}\label{eq:Green}
 P_\theta(z(\omega))
 =\int_{\partial\Om}
       \bigl(u\partial_\nu v-v\partial_\nu u\bigr)\dd S.
\end{equation}
Since $u=0$ on $\partial\Om$, the term $u\partial_\nu v$ vanishes. The Cauchy--Schwarz inequality gives
\[
 |P_\theta(z(\omega))|\le |\partial\Om|^{1/2}
 e^{R_\Om|\Ima z(\omega)|}d(\omega),
\]
and hence
\begin{equation}\label{eq:curvedata}
 |P_\theta(z(\omega))|
 \le C(1+|z(\omega)|)e^{R_\Om|\Ima z(\omega)|}d(\omega).
\end{equation}

\begin{lemma}\label{lem:interpolation}
Suppose $P_\theta$ satisfies~\eqref{eq:entiregrowth} and~\eqref{eq:curvedata}, with $\int_0^K d(\omega)^2\dd\omega=\varepsilon_K^2$. If $k=\sqrt{bK}\ge k_0=\max\{1,\sqrt a\}$, then we have the following estimate which holds for every $r\ge k_0$:
\begin{equation}\label{eq:lowfourier}
 \sup_{|z|\le r}|P_\theta(z)|
 \le C\varepsilon_K e^{C r\log(2+r/k)}+CM e^{-2r},
\end{equation}
where $C$  is a generic constant independent of $\theta,K,r,M$, and $\varepsilon_K$.
\end{lemma}

\begin{proof}

Fix $r\ge k_0$ and set
\begin{equation}\label{eq:Tnodes}
 T=\min\{K,r^2/b\},\qquad t=\sqrt{bT}=\min\{k,r\}.
\end{equation}
Then $1\le t\le r$ and $a\le t^2$. For an integer $N\ge2$, consider
\[
 J_j=\left[\frac T2+\frac{(j-1)T}{2N},\,
                \frac T2+\frac{(j-1)T}{2N}+\frac T{4N}\right],
 \qquad 1\le j\le N.
\]
Each interval has length $T/(4N)$ and lies in $[T/2,T]$. Since
\[
 \frac1{|J_j|}\int_{J_j}d(\omega)^2\dd\omega
 \le\frac{4N}{T}\varepsilon_K^2,
\]
we can choose $\omega_j\in J_j$ for which
\begin{equation}\label{eq:nodedata}
 d(\omega_j)\le2\sqrt{N/T}\,\varepsilon_K.
\end{equation}
The function $d$ is continuous by the resolvent representation and the trace theorem, so the averaging choice is valid also at zero discrepancy. The same nodes are used for all $\theta$. Then we put $z_j=\sqrt{-a-\ii b\omega_j}$ on the branch specified above. We have
\[
|z_j|=(a^2+b^2\omega_j^2)^{1/4}\le2^{1/4}t\le2t.
\]
For $j>m$, the gap between the intervals gives
\[
 \omega_j-\omega_m
 \ge\frac{T}{4N}\bigl(2(j-m)-1\bigr)
 \ge\frac{T}{4N}(j-m).
\]
Using $z_j^2-z_m^2=-\ii b(\omega_j-\omega_m)$ and
$|z_j+z_m|\le4t$, we obtain
\begin{equation}\label{eq:nodegeom}
 |z_j|\le2t,\quad{\rm and}\quad
 |z_j-z_m|
 =\frac{b|\omega_j-\omega_m|}{|z_j+z_m|}
 \ge\frac{t}{16N}|j-m|\quad{\rm for}~j\ne m.
\end{equation}
In particular these complex nodes are distinct. The trace estimate and~\eqref{eq:nodedata} imply
\begin{equation}\label{eq:nodeP}
 \max_j|P_\theta(z_j)|
 \le C(1+t)e^{2R_\Om t}\sqrt{N/T}\,\varepsilon_K
 \le C e^{2R_\Om r}\sqrt N\,\varepsilon_K.
\end{equation}
For the last inequality we have used $\sqrt T=t/\sqrt b$, $t\ge1$, and $(1+t)/t\le2$.

Next, we construct the interpolation polynomial, which is
a polynomial of degree at most $N-1$ whose value at $z_j$ is $P_\theta(z_j)$ for every $j$. Define the Lagrange cardinal polynomials
\[
 \ell_j(z)=\prod_{\substack{1\le m\le N\\m\ne j}}
                   \frac{z-z_m}{z_j-z_m},
 \qquad
 \Pi_{N-1}(z)=\sum_{j=1}^N P_\theta(z_j)\ell_j(z),
\]
where
\[
 \ell_j(z_m)=\delta_{jm}
 =\begin{cases}1,&j=m,\\0,&j\ne m.\end{cases}
\]
Thus, we have $\Pi_{N-1}(z_j)=P_\theta(z_j)$. The polynomial is unique, since the difference of two such polynomials has $N$ distinct zeros and degree at most $N-1$. The same argument shows that every polynomial $q$ of degree at most $N-1$ has the representation
$q(z)=\sum_jq(z_j)\ell_j(z)$.
Consequently $\ell_1,\ldots,\ell_N$ form a basis of this $N$-dimensional polynomial space. This basis depends on the nodes, and orthogonality is neither assumed nor used. For details on the Lagrange formulas we refer to \cite[\S3.3(i), eqs.~(3.3.1)--(3.3.3)]{DLMF}.

If $|z|\le r$, then $|z-z_m|\le r+2t$. Also, \eqref{eq:nodegeom} yields
\[
 \prod_{m\ne j}|z_j-z_m|
 \ge\left(\frac{t}{16N}\right)^{N-1}
       (j-1)!(N-j)!.
\]
Indeed, the factors $|j-m|$ are $1,\ldots,j-1$ on one side and $1,\ldots,N-j$ on the other. Then it follows that
\begin{align*}
 \sum_{j=1}^N|\ell_j(z)|
 &\le\left(\frac{16N(r+2t)}{t}\right)^{N-1}
       \sum_{j=1}^N\frac1{(j-1)!(N-j)!}\\
 &=\left(\frac{16N(r+2t)}{t}\right)^{N-1}
       \frac{2^{N-1}}{(N-1)!}\\
 &\le\left(\frac{64e(r+2t)}{t}\right)^{N-1}
 \le\bigl(C_{\mathrm L}(1+r/t)\bigr)^N,
 \qquad C_{\mathrm L}=128e.
\end{align*}
The equality uses the binomial identity
$\sum_{j=1}^N\binom{N-1}{j-1}=2^{N-1}$.
For the next inequality, we have used $m!\ge(m/e)^m$ with $m=N-1$ and $N/(N-1)\le2$.
The factorial bound follows from
\[
 \log(m!)=\sum_{h=1}^m\log h
 \ge\int_1^m\log x\dd x
 =m\log m-m+1\ge m\log m-m.
\]
Combining the basis bound with~\eqref{eq:nodeP}, we obtain
\begin{equation}\label{eq:polybound}
 |\Pi_{N-1}(z)|
 \le C e^{2R_\Om r}\sqrt N\,
       \bigl(C_{\mathrm L}(1+r/t)\bigr)^N\varepsilon_K,
 \qquad |z|\le r.
\end{equation}

Let $W_N(w)=\prod_{j=1}^N(w-z_j)$ be the nodal polynomial. On the counterclockwise circle $|\zeta|=8r$, all nodes and every point $|z|\le r$ are strictly inside the contour. 
Then the remainder formula is
\begin{equation}\label{eq:remainder}
 P_\theta(z)-\Pi_{N-1}(z)
 =\frac1{2\pi\ii}\int_{|\zeta|=8r}
 \frac{P_\theta(\zeta)}{\zeta-z}
 \prod_{j=1}^N\frac{z-z_j}{\zeta-z_j}\dd\zeta,
\end{equation}
The derivation of the formula is as follows.
First assume $z\notin\{z_1,\ldots,z_N\}$. As a function of $\zeta$, the integrand
\[
 H_z(\zeta)=\frac{P_\theta(\zeta)W_N(z)}
                  {(\zeta-z)W_N(\zeta)}
\]
has no possible singularities other than simple poles at $z$ and the nodes. Its residues are
\begin{align*}
 \operatorname{Res}_{\zeta=z}H_z(\zeta)&=P_\theta(z),\\
 \operatorname{Res}_{\zeta=z_j}H_z(\zeta)
 &=\frac{P_\theta(z_j)W_N(z)}{(z_j-z)W_N'(z_j)}
 =-P_\theta(z_j)\ell_j(z).
\end{align*}
The last equality follows from
\[
 W_N'(z_j)=\prod_{m\ne j}(z_j-z_m),\qquad
 \frac{W_N(z)}{(z-z_j)W_N'(z_j)}=\ell_j(z).
\]
The residue theorem therefore gives~\eqref{eq:remainder}, which extends to $z=z_j$ by continuity with both sides equal to zero. This is also the contour remainder formula in~\cite[\S3.3(i), eq.~(3.3.6)]{DLMF}.

On the contour and for $|z|\le r$, the bounds $t\le r$ and $|z_j|\le2t$ give
\[
 |\zeta-z|\ge7r,\qquad |\zeta-z_j|\ge6r,
 \qquad |z-z_j|\le3r.
\]
Moreover, we have $|P_\theta(\zeta)|\le|\Om|^{1/2}Me^{8R_\Om r}$ by~\eqref{eq:entiregrowth}. Since the circle has length $16\pi r$, formula~\eqref{eq:remainder} implies
\[
 |P_\theta(z)-\Pi_{N-1}(z)|
 \le\frac{16\pi r}{2\pi\,7r}
       |\Om|^{1/2}Me^{8R_\Om r}\left(\frac{3r}{6r}\right)^N
 =\frac87|\Om|^{1/2}Me^{8R_\Om r}2^{-N}.
\]

Now we fix $A_*>0$ with $A_*\log2\ge8R_\Om+2$, and choose
\[
 N=\lceil A_*r\rceil+2.
\]
Here $\lceil x\rceil$ denotes the least integer greater than or equal to $x$. Then one has $e^{8R_\Om r}2^{-N}\le\tfrac14 e^{-2r}$, and thus the remainder is at most $CM e^{-2r}$. Also, since $r\ge1$, we have
$N\le(A_*+3)r$. The logarithm of the amplification factor in~\eqref{eq:polybound}, apart from its fixed leading constant, is
\[
 2R_\Om r+\tfrac12\log N
       +N\log\bigl(C_{\mathrm L}(1+r/t)\bigr)
 \le Cr\log(2+r/t).
\]
Here $A_*$ and $C_{\mathrm L}$ are fixed independently of $r,k$, and $N$, and $r/t\ge1$. Thus
\[
 |P_\theta(z)|\le C\varepsilon_K e^{Cr\log(2+r/t)}+CM e^{-2r},
 \qquad |z|\le r.
\]
Finally, for $r\ge k$, we let $t=k$, and if $r<k$, we set $t=r$ and then
\[
 \log(2+r/t)=\log3
 \le\frac{\log3}{\log2}\log(2+r/k).
\]
The constants are uniform in $\theta$. Taking the indicated suprema completes the proof.
\end{proof}

In the following theorem, we establish an increasing stability for the inverse source problem. The stability shows that as the upper bound of the frequency band or regularity of the source increases, the stability will also improve.
\begin{theorem}\label{thm:constant}
Let $p\in\mathcal B_M$, $s>0$, $k_0=\max\{1,\sqrt a\}$, and $k=\sqrt{bK}\ge k_0$. For $M>0$ and $\varepsilon_K>0$, set $E_K=\log(2+M/\varepsilon_K)$. One has the stability estimate
\begin{equation}\label{eq:mainconst}
 \norm p_{L^2(\R^d)}
 \le CM\left(1+\frac{E_K}{\log(2+E_K/k)}\right)^{-s}.
\end{equation}
The constant $C$ is independent of $K,M$, and $\varepsilon_K$. 
\end{theorem}
\begin{proof}
Use the unitary Fourier transform $\widehat p(\xi)=(2\pi)^{-d/2}\int_{\R^d}p(x)e^{-\ii\xi\cdot x}\dd x$. For $\xi\ne0$, $\widehat p(\xi)=(2\pi)^{-d/2}P_{\xi/|\xi|}(|\xi|)$. For every $r\ge k_0$, Plancherel's theorem,~\eqref{eq:fourierprior}, and Lemma~\ref{lem:interpolation} give
\begin{align}\label{eq:fouriersplit}
 \norm p_{L^2}^2&=\int_{|\xi|\leq r}|\widehat p(\xi)|^2\dd\xi + \int_{|\xi|> r}|\widehat p(\xi)|^2\dd\xi\notag\\
 &\le C r^d\varepsilon_K^2e^{C r\log(2+r/k)}
     +CM^2r^de^{-4r}
     +\int_{|\xi|>r}|\widehat p(\xi)|^2\dd\xi\\
 &\le C\varepsilon_K^2e^{C r\log(2+r/k)}+CM^2r^{-2s}.
 \nonumber
\end{align}

Here the factor $r^d$ in the data term is absorbed by increasing the constant in the exponential, uniformly in $k$, since $r\ge1$ and $\log(2+r/k)\ge\log2$. Also $r^de^{-4r}\le C_s r^{-2s}$, and the Sobolev prior controls the Fourier tail by $M^2r^{-2s}$.

Taking square roots in~\eqref{eq:fouriersplit} and changing the fixed constant yields
\[
 \norm p_{L^2}\le C\varepsilon_K e^{Cr\log(2+r/k)}+CMr^{-s},
 \qquad r\ge k_0.
\]
Write $E=E_K$ and $Q=E/\log(2+E/k)$. Choose a fixed $c>0$ with $c\le\log2$ and $Cc\le1/4$, and put $Q_*=\max\{k_0/c,2\}$. If $Q\le Q_*$, the prior immediately gives
\[
 \norm p_{L^2}\le M\le (1+Q_*)^sM(1+Q)^{-s}.
\]
If $Q>Q_*$, take $r=cQ$. Then $r\ge k_0$, $r\le E$, and $E\ge Q\log2>\log4$. Consequently
\[
 Cr\log(2+r/k)\le CcE\le E/4,
 \qquad \varepsilon_K/M=(e^E-2)^{-1}\le2e^{-E}.
\]
The preceding unsquared estimate becomes
\[
 \norm p_{L^2}\le CM e^{-3E/4}+CMQ^{-s}
 \le C_sM(1+Q)^{-s}.
\]
The last step uses $Q\le E/\log2$ and $Q>Q_*\ge2$. All thresholds and constants are independent of $k\ge k_0$, $M$, and $\varepsilon_K$. This proves~\eqref{eq:mainconst} in each case. The proof is complete.
\end{proof}

\section{Numerical reconstruction}\label{sec:numerics}

We describe a reconstruction scheme for the constant-coefficient problem \eqref{eq:constant} in $\Om=B_1(0)\subset\R^3$, with $a=b=1$ and homogeneous Dirichlet boundary data. The source is real-valued and the measured quantity is
$q(\omega)=\partial_\nu u(\omega)$ on $\partial\Om$.
We follow the frequency-by-frequency Landweber strategy of Bao, Lin, and Triki~\cite{BLT2011} and Bao et~al.~\cite{BLRX2015}; see also the
survey~\cite{BLLT2015} and the later application in~\cite[Section~5]{IL}. For the present diffusion problem, we use repeated sweeps and a dual-residual discrepancy rule.
Here the spherical geometry permits a spherical-harmonic representation.

\subsection{Forward calculation and reference data}
Let $Y_{\ell m}$, $-\ell\le m\le\ell$, be real orthonormal spherical harmonics
on $\mathbb S^2$, and write $p(r,\theta)=\sum_{\ell,m}p_{\ell m}(r)Y_{\ell m}(\theta)$.
For the numerical examples, assume $\supp p\subset\overline B_{R_s}$, where
$R_s=0.8$. If $\gamma_\omega=\sqrt{a+\ii b\omega}$ with positive real part,
Green's identity gives
\begin{equation}\label{eq:num-green}
	q_{\ell m}(\omega)
	=-\int_0^{R_s}p_{\ell m}(r)
	\frac{i_\ell(\gamma_\omega r)}{i_\ell(\gamma_\omega)}r^2\dd r,
\end{equation}
where $i_\ell$ is the modified spherical Bessel function. The minus sign corresponds to the outward normal derivative. Reference data are generated from the explicitly defined smooth source by angular quadrature and high-order radial
quadrature in~\eqref{eq:num-green}. Scaled Bessel functions are used to evaluate the quotient at large frequencies.

For inversion, let $\phi_{\ell n}$ be the radial $L^2((0,R_s),r^2\dd r)$ orthonormalization of
$\chi(r)j_\ell(z_{\ell n}r/R_s)$, $1\le n\le N_r$, where $z_{\ell n}$ is the $n$th positive zero of the spherical Bessel function $j_\ell$, and $\chi(r)=\exp(-r^2/(R_s^2-r^2))$ for $r<R_s$, extended by zero.
Set
\begin{equation}\label{eq:num-source}
	p_h(r,\theta)=\sum_{\ell=0}^{L}\sum_{m=-\ell}^{\ell}
	\sum_{n=1}^{N_r}c_{\ell mn}\phi_{\ell n}(r)Y_{\ell m}(\theta),
	\qquad c_{\ell mn}\in\R.
\end{equation}
This is a numerical approximation space of dimension $(L+1)^2N_r$; its functions have smooth, compactly supported zero extensions. The smooth reference source is not replaced by its projection onto this space.

The forward matrix is assembled independently by a spectral Galerkin method.
Let $R_{\ell t}(r)$ be the Dirichlet radial eigenfunction
$j_\ell(z_{\ell t}r)$ normalized by $\int_0^1R_{\ell t}(r)^2r^2\dd r=1$. For frequency $\omega_j$, its row acting
on $c_{\ell m}=(c_{\ell m1},\ldots,c_{\ell mN_r})^T$ is
\begin{equation}\label{eq:num-forward}
	(A_{\ell,j})_n=
	\sum_{t=1}^{N_f}\frac{R'_{\ell t}(1)}{z_{\ell t}^2+a+\ii b\omega_j}
	\int_0^{R_s}R_{\ell t}(r)\phi_{\ell n}(r)r^2\dd r.
\end{equation}
Thus $q_{\ell m}(\omega_j)\approx A_{\ell,j}c_{\ell m}$; the signed derivative $R'_{\ell t}(1)$ is retained. Formula~\eqref{eq:num-forward} follows by expanding $u$ in the Dirichlet eigenfunctions and dividing each source coefficient by $z_{\ell t}^2+a+\ii b\omega_j$. Increasing $N_f$ and the quadrature orders, and comparing with~\eqref{eq:num-green} on trial functions, checks the forward accuracy without using reconstruction errors.

\subsection{Multi-frequency reconstruction}
Choose nested midpoint grids $\omega_j=(j-\tfrac12)\Delta\omega$, $j=1,\ldots,J_K$, with $K=J_K\Delta\omega$ and $w_j=\Delta\omega$.
Retain observation harmonics through degree $L_o\ge L$. The data norm is
\begin{equation}\label{eq:num-data-norm}
	\norm d_{W,K}^2=\sum_{j=1}^{J_K}w_j
	\sum_{\ell=0}^{L_o}\sum_{m=-\ell}^{\ell}|d_{j,\ell m}|^2,
\end{equation}
which approximates the frequency-integrated boundary norm in
\eqref{eq:epsilon}. Set $d^\delta=q^{\mathrm{ref}}+\xi$, where
$q^{\mathrm{ref}}$ denotes the computed reference data. Let $A_j$ collect the rows in~\eqref{eq:num-forward}, with zero predicted components for $L<\ell\le L_o$. For real coefficients, define $B_j=(\Rea A_j;\Ima A_j)$ and $y_j=(\Rea d_j^\delta;\Ima d_j^\delta)$ by vertical concatenation. Then the real adjoint is
$B_j^Ty_j=\Rea(A_j^*d_j^\delta)$.

To keep the iteration settings fixed as $K$ varies, construct one metric from the known operator on the largest prescribed band $K_*$. Let $H$ be the Gram matrix of the whole-space $H^2$ inner product on the basis in
\eqref{eq:num-source}, so that $c^THc=\norm{p_h}_{L^2}^2+2\norm{\nabla p_h}_{L^2}^2+\norm{\Delta p_h}_{L^2}^2$.
Set
\begin{equation}\label{eq:num-metric}
	G_* =\sum_{j=1}^{J_{K_*}}w_jB_j^TB_j,\qquad
	\beta=\lambda_{\max}(H^{-1/2}G_*H^{-1/2}),\qquad
	\mathsf{M}=G_*+10^{-2}\beta H.
\end{equation}
The positive definite matrix $\mathsf{M}$ is computed without observation data and reused for all bands. It defines an iteration metric, not a bound on $\norm{p_h}_{H^2}$. Define, once on the largest grid,
\begin{equation}\label{eq:num-steps}
	\nu_j=w_j\norm{B_j\mathsf{M}^{-1/2}}_2^2,\qquad
	\nu_* =\max_{1\le j\le J_{K_*}}\nu_j,\qquad
	\mu_j=\frac{0.2}{\nu_*}.
\end{equation}
The implemented step-floor ratio is one, so the step is fixed across frequency blocks. At each $K$, start from $c=0$ and sweep through $j=1,\ldots,J_K$ in ascending order. At frequency $j$, perform two updates
\begin{equation}\label{eq:num-kl}
	c\ \leftarrow\ c+\mu_jw_j\mathsf{M}^{-1}B_j^T(y_j-B_jc).
\end{equation}
The choice~\eqref{eq:num-steps} ensures $\mu_j\nu_j\le0.2<2$ for every block.
After complete sweeps, stop at the first coefficient vector satisfying
\begin{equation}\label{eq:num-stop}
	\|r_K(c)\|_2\le1.05(\delta_K+\epsilon_f),\qquad
	\|P_Kr_K(c)\|_2\le1.05(\kappa_K\delta_K+\epsilon_f),
\end{equation}
where $r_K(c)$ is the weighted realification of $A_Kc-d^\delta$,
$(A_Kc)_j=A_jc$, and $\delta_K=\|\xi_K\|_{W,K}$.
The data-independent orthogonal projection $P_K$ retains discrete Legendre degrees $0,\ldots,3$ on frequency segments of width $4$ for harmonics $\ell\le8$, separately for real and imaginary parts; a segment with $n$ nodes retains
$\min\{4,n\}$ orthonormal polynomial coordinates.
With $D_K=2J_K(L_o+1)^2$ and $d_K=\operatorname{rank}P_K$,
\[
\kappa_K^2=F^{-1}_{d_K/2,(D_K-d_K)/2}(0.99),
\]
where $F_{a,b}$ is the beta-distribution cumulative distribution function.
This calibrates the projected Gaussian noise direction per case, rather than a reconstruction confidence bound. Raw data drive all updates.

The numerical allowance is $\epsilon_f=2(\epsilon_{\rm ref}+R_{\rm num}\epsilon_{\rm op})$,
$R_{\rm num}=1$, using independent reference and operator refinement differences. It is an empirical estimate, subject to a coefficient-radius check, not a certified continuous error bound. No uniform source-space mismatch bound was declared; the known-source projection defect is reported separately and does not determine stopping.
All observation channels enter the full residual, which need not decrease monotonically over successive sweeps. If both targets are not reached within $2\times10^7$ complete sweeps, the run is reported as unsuccessful.
Finite-dimensional truncation and discrepancy stopping provide regularization. For efficiency, complete sweeps may be composed into affine maps and evaluated using matrix powers, retaining the first verified dual-residual hit on the Kaczmarz–Landweber path in exact arithmetic. Floating-point checks are not interval-arithmetic certification.

\subsection{Bandwidth comparison and accuracy checks}
The baseline discretization uses $L=N_r=8$, $L_o=12$, $N_f=160$ and $\Delta\omega=0.125$, giving 648 source coefficients. We take $K=m^2$, $m=2,\ldots,16$, and keep the source space, iteration metric and step fixed throughout the bandwidth comparison. The connected-bump, ellipsoidal and toroidal test sources are smooth, compactly supported in $\overline B_{0.8}$ and normalized to unit $L^2$ norm. Reference data are generated by direct integration of the prescribed sources.

\begin{figure}[htbp]
	\centering
	\includegraphics[width=\textwidth]{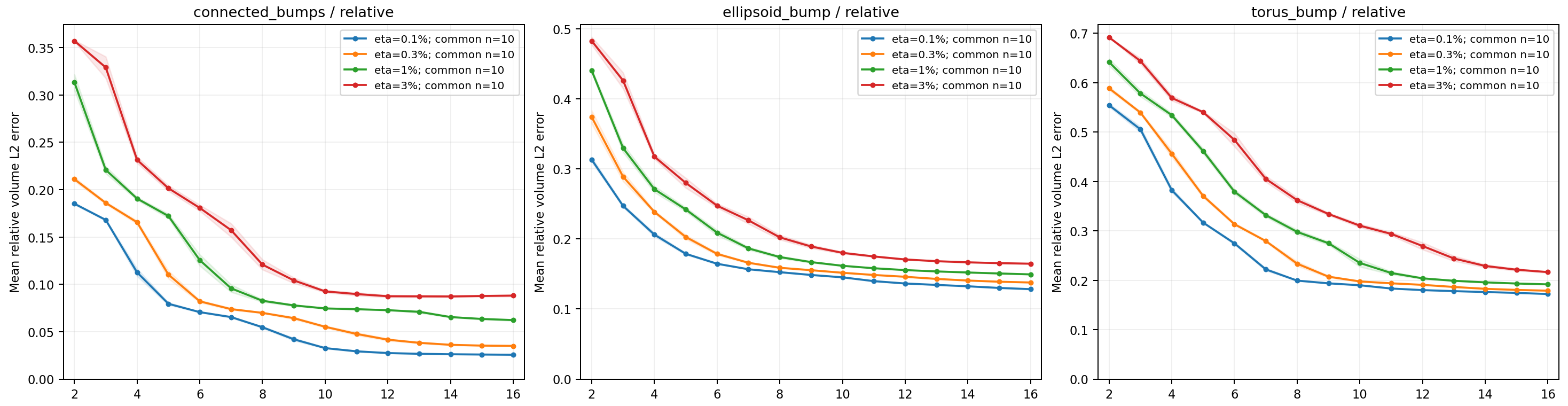}
	\caption{Mean relative volume $L^2$ error versus $\sqrt K$ for three
		sources and four relative noise levels.}
	\label{fig:num-main}
\end{figure}

Noise is generated from complex Gaussian directions and normalized in each band to the relative level $\eta\in\{10^{-3},3\times10^{-3},10^{-2},3\times10^{-2}\}$. We use ten paired trials across all bands and report mean relative volume $L^2$ errors evaluated against the full true sources. Shading denotes pointwise 95\% Student-$t$ confidence intervals for the mean error.
Figure~\ref{fig:num-main} shows that increasing the observation bandwidth leads to an overall reduction in the mean reconstruction errors for all three sources at the tested noise levels. At 1\% noise, as $K$ increases from $4$ to $256$, the mean errors for the connected, ellipsoidal and toroidal sources decrease from 31.34\%, 44.00\% and 64.15\% to 6.22\%, 14.93\% and 19.16\%, respectively.

\begin{figure}[H]
	\centering
	\includegraphics[width=\linewidth]{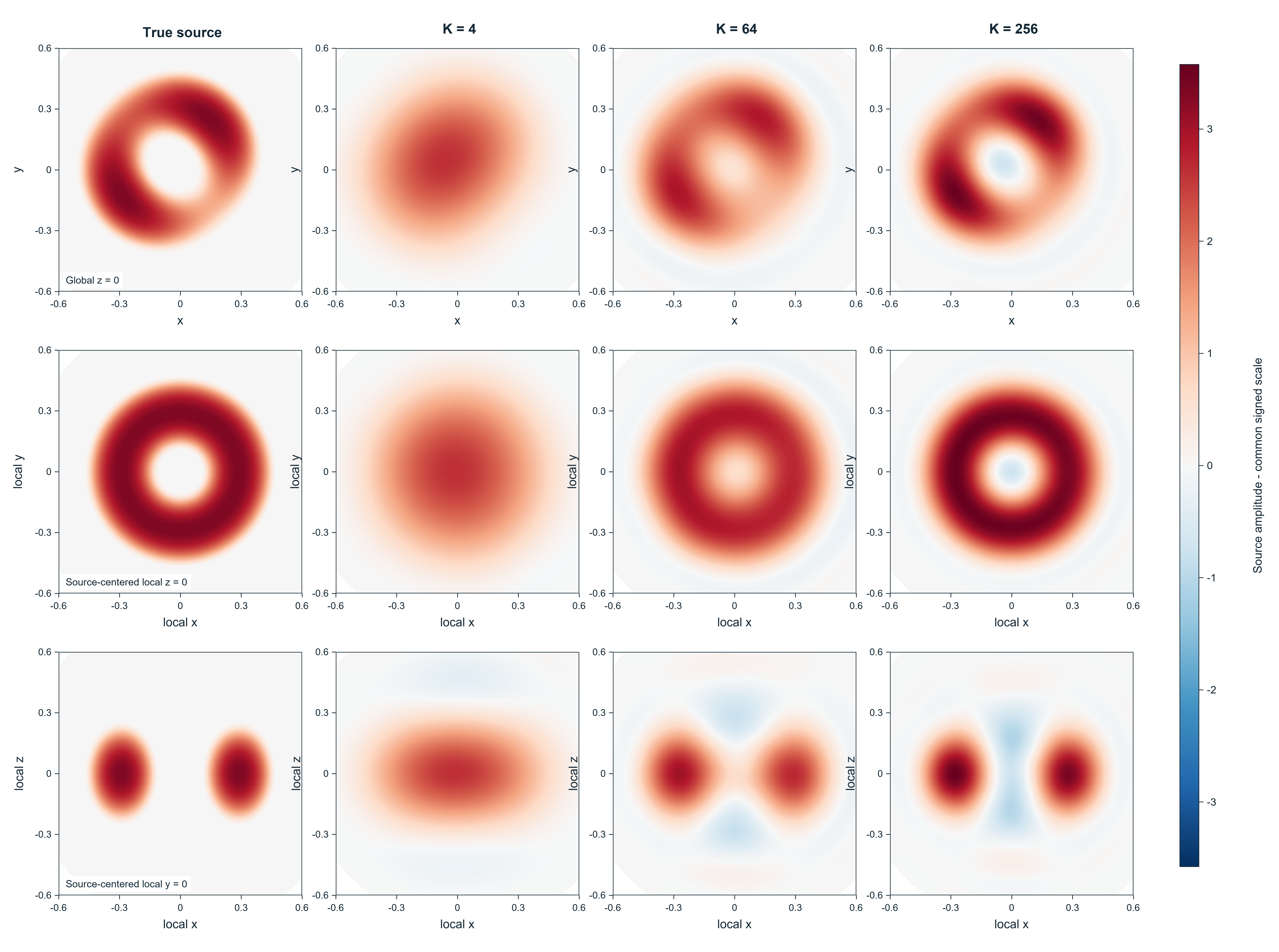}
	\caption{True toroidal source and reconstructions at $1\%$ relative noise.}
	\label{fig:num-torus}
\end{figure}

Figure~\ref{fig:num-torus} compares the true toroidal source with reconstructions from a single trial at $K=4,64,256$. The reconstructions become clearer as the observation bandwidth increases.

The error curves and reconstruction images show that increasing the observation bandwidth reduces reconstruction errors overall and improves the resolution of source structures, consistent with the bandwidth dependence established in Theorem~\ref{thm:constant}.

\section{Conclusion}\label{sec:conclusion}

We studied the recovery of a frequency-independent source in the time-harmonic diffusion equation from finite-band boundary flux measurements. Under the prescribed a priori source bounds, heat observability and quantitative analytic continuation yield uniqueness and a conditional double-logarithmic stability estimate for variable coefficients. In the constant-coefficient case, interpolation of the Fourier--Laplace transform provides an increasing stability estimate with explicit dependence on the observation bandwidth and the Sobolev regularity of the source.

For the three-dimensional constant-coefficient problem in the unit ball, source reconstruction is performed in a fixed approximation space using a multi-frequency Landweber iteration with a dual-residual discrepancy stopping rule. Numerical results for three smooth, compactly supported sources show that increasing the observation bandwidth reduces reconstruction errors overall and improves the resolution of source structures under the tested relative-noise protocol. These results are consistent with the bandwidth dependence established in Theorem~\ref{thm:constant}.

\section*{Acknowledgements}
 X. Xu is partially supported by the National Key Research and Development Program of China (No. 2024YFA1012303), National Natural Science Foundation of China (No. 12525112), and the Open Research Project of Innovation Center of Yangtze River Delta, Zhejiang University. Y. Zhao is supported by National Natural Science Foundation of China (No. 12371421).

\section*{Conflict of Interest Statement}

The authors declare that they have no conflict of interest.

\section*{Data Availability Statement}
The numerical experiments use synthetic data generated using the models and parameters described in this article.
The data supporting the numerical results are available from the corresponding author upon reasonable request.

\end{document}